\documentclass[conference]{IEEEtran}
\IEEEoverridecommandlockouts
\usepackage{cite}
\usepackage{amsmath,amssymb,amsfonts}
\usepackage{algorithmic}
\usepackage{graphicx}
\usepackage{textcomp}
\usepackage{xcolor}
\def\BibTeX{{\rm B\kern-.05em{\sc i\kern-.025em b}\kern-.08em
		T\kern-.1667em\lower.7ex\hbox{E}\kern-.125emX}}
\usepackage{cite}
\usepackage{multirow}
\usepackage{bigstrut}
\usepackage{comment}
\usepackage[normalem]{ulem}
\usepackage{amsmath,amssymb,amsfonts}
\usepackage{algorithmic}
\usepackage{algorithm}
\usepackage{graphicx}
\usepackage{textcomp}
\usepackage{xcolor}

\newtheorem{theo}{Theorem}[section]

\newtheorem{rem}{Remark} 

\usepackage{dsfont}

\DeclareMathOperator*{\argmax}{arg\,max}
\DeclareMathOperator*{\argmin}{arg\,min}
\newcommand{\proof}     {\paragraph*{Proof}}
\newcommand{\carre}     {\hfill$\Box$}

\begin{document}
	
	\title{Exponentially Weighted Algorithm for Online Network Resource Allocation with Long-Term Constraints  
		\thanks{This work was supported by a grant from the Natural Sciences and Engineering Research Council of Canada and Ericsson Canada. }
	}
	
	\author{\IEEEauthorblockN{ Ahmed Sid-Ali\IEEEauthorrefmark{1}, Ioannis Lambadaris\IEEEauthorrefmark{2}, Yiqiang Q. Zhao\IEEEauthorrefmark{1}, Gennady Shaikhet\IEEEauthorrefmark{1}, and Amirhossein Asgharnia\IEEEauthorrefmark{2}}
		\IEEEauthorblockA{\IEEEauthorrefmark{1}School of Mathematics and Statistics
			\\}
		\IEEEauthorblockA{\IEEEauthorrefmark{2}Department of Systems and Computer Engineering\\
			Carleton University, Ottawa, Ontario\\
			Emails: Ahmed.Sidali@carleton.ca; ioannis@sce.carleton.ca; zhao@math.carleton.ca; gennady@math.carleton.ca;\\ 
			and  amirhosseinasgharnia@sce.carleton.ca}
		
	}

	\maketitle
	
	\begin{abstract}
		This paper studies an online optimal resource reservation problem in communication networks with job transfers where the goal is to minimize the reservation cost while maintaining the blocking cost under a certain budget limit. To tackle this problem, we propose a novel algorithm based on a randomized exponentially weighted method that encompasses long-term constraints. We then analyze the performance of our algorithm by establishing an upper bound for the associated \textit{regret} and for the cumulative constraint violations. Finally, we present numerical experiments where we compare the performance of our algorithm with those of reinforcement learning where we show that our algorithm surpasses it.
	\end{abstract}
	
	
	\section{Introduction}
	Online optimization is a machine learning framework where a decision maker sequentially chooses a sequence of decision variables over time to minimize the sum of a sequence of (convex) loss functions. The decision maker does not have full access to the data at once but receives it incrementally. Thus, the decisions are made sequentially in an incomplete information environment. Moreover, no assumptions on the statistics of data sources are made. This radically differs from classical statistical approaches such as Bayesian decision theory or Markov decision processes. Therefore, only the observed data sequence's empirical properties matter, allowing us to address the dynamic variability of traffic requests in modern communication networks. Online optimization has further applications in a wide range of fields; see, e.g. \cite{b12, b13, b14} and the references therein for an overview. Moreover, given that the decision-maker can only access limited/partial information, globally optimal solutions are generally not realizable. Instead, one searches for algorithms that perform relatively well compared to the overall {\it ideal} best solution in hindsight which has full access to the data. This performance metric is referred to as \textit{regret} in the literature. In particular, if an algorithm incurs regret that increases sub-linearly with time, then one says that it achieves \textit{no regret}. 
	
	One important application of online optimization, as in the problem presented in our work, is in resource reservations that arise in situations where resources must be allocated in advance to meet future unknown demands. In this paper, we tackle a specific online resource reservation problem in a simple communication network topology. The network is composed of $N$ linked servers and the network administrator reserves resources at each server to meet future job requests. The specificity of this system is that it allows to transfer of jobs from one server to another after receiving the job requests to meet the demands better. This couples the servers by creating dependencies between them. The reservation and transfer steps come with a cost that is proportional to the amount of resources involved. Moreover, a violation cost is incurred if the job requests are not fully satisfied. The goal is then to minimize the reservation cost while maintaining the transfer and violation costs under a given budget threshold which leads to an online combinatorial optimization problem with long-term constraints. To tackle it, we propose a randomized control policy that minimizes the cumulative reservation cost while maintaining the long-term average of the cumulative expected violation and transfer costs under the budget threshold. In particular, we propose a novel exponentially weighted algorithm that incorporates long-term constraints. We then derive an explicit upper bound for the incurred \textit{regret} and an upper bound for the cumulative constraints violations. 
	
	The rest of the paper is organized as follows: in Section \ref{model-sec} we introduce the problem; in Section \ref{rand-formulation-sec}, we formalize the problem as a constrained online optimization problem on the simplex of probability distributions over the space of reservations; section \ref{algo-sec} contains our proposed exponentially weighted algorithm with constraints; then, in Section \ref{per-sec}, we present an upper bound for the regret together with the cumulative constraint violations upper bound, and finally, we present in Section \ref{num-sec} some numerical results where we compare the performance of our algorithm with a tailored reinforcement learning approach.  
	
	\section{Online resource reservation in communication networks}  
	\label{model-sec} 
	Consider a network of $N$ servers connected by communication links. The system provides access to computing resources for clients. For simplicity, suppose that there is a single type of resource (e.g. memory, CPU, etc.). Denote by $m_n$  the total resources available at the $n$-th server. Let $\mathcal{R}_n=\{1,\ldots,m_n\}$ be the set of possible reservations at the $n$-th server, and $\mathcal{R}=\prod_{n=1}^N\mathcal{R}_n$ be the set of possible reservations in the entire network. The system operates in discrete time where, at each time slot $t=1,2,\ldots$, the following processes take place chronologically: 
	\begin{itemize} 
		\item \textbf{Resource reservation}: the network administrator selects the resources $A^t=(A^t_1,\ldots,A^t_N)\in\mathcal{R}$ to make available at each server.    
		\item \textbf{Job requests}: the network receives job requests $B^t=(B^t_1,\ldots,B^t_N)\in\mathcal{R}$ from its clients to each of its servers. 
		\item \textbf{Job transfer}: the network administrator can shift jobs between the servers to best accommodate the demands.  
	\end{itemize} 
	Let $\delta^t_{n,m}$ be the number of jobs transferred from server $n$ to server $m$ at time slot $t$ after receiving the job requests $B^t$. These quantities depend on both $A^t$ and $B^t$, however, this dependency is suppressed to keep the notation simple. Moreover, suppose that the following costs are associated with each of the three processes:
	\paragraph*{\textbf{-Reservation cost}} given by the function
	\begin{align*}
		C(A^t)=\sum_{n=1}^N  f^R_n(A^t_n),
	\end{align*} 
	
	\paragraph*{\textbf{-Violation cost}} incurred  when the job requests cannot be satisfied  
	\begin{align*} 
		C_V(A^t,B^t)=\sum_{n=1}^N f^V_n\big(B_n^t-A_n^t-\sum_{m=1}^N\delta^t_{n,m}\big),
	\end{align*} 
	
	\paragraph*{\textbf{-Transfer cost}} incurred by the transfer of jobs between the different servers 
	\begin{align*}
		&C_T(A^t,B^t)=\sum_{n=1}^N\sum_{m\neq n} f^T_n\big(\delta^t_{n,m}\big), 
	\end{align*} 
	where $f^R_n,f^V_n,f^T_n$, for $1\leq n\leq N$, are some positive functions. Therefore, at each time slot $t$, and after receiving the job requests $B^t$, the job transfer coefficients $\{\delta^t_{n,m},1\leq n \leq N\}$ are solutions to the following offline minimization problem: 
	\begin{equation}  
		\begin{split}
			\left\{
			\begin{tabular}{l}
				$\min \sum\limits_{n=1}^N\big(\sum\limits_{m\neq n}f^T_n(\delta^t_{n,m})+f^V_n\big(B_n^t-A_n^t-\sum\limits_{m\neq n}\delta^t_{n,m}\big)\big),$\\
				s.t. $\delta^t_{n,m}\leq\min\big\{ (B^t_n-A_n^t)^+,(A^t_{m}-B^t_{m})^+ \big\},\forall  n\neq m$. 
			\end{tabular}
			\right.
			\label{delta-opt-prob}
		\end{split}
	\end{equation} 
	One then aims to minimize the reservation cost $C(A^t)$ while maintaining the sum  $C_0(A^t,B^t):=C_V(A^t, B^t)+C_T(A^t, B^t)$ of the violation and transport costs under a given threshold $v>0$. Nevertheless, since it is an online problem, its optimal solution is out of reach. In particular, one cannot guarantee that the constraint is satisfied. Instead, one aims to solve the following online optimization problem with a long-term constraint over time horizons $T>0$:   
	\begin{equation} 
		\begin{split}
			\left\{
			\begin{tabular}{l}
				$\min\limits_{\{A^t\}_{t=1}^T}\sum\limits_{t=1}^T C(A^t)$,\\
				s.t.:  $ \frac{1}{T} \sum\limits_{t=1}^T C_0(A^t,B^t)\leq v$, 
			\end{tabular}
			\right.
			\label{opt-pb-on}
		\end{split}
	\end{equation}
	More precisely, one searches for an online control policy that uses the cumulative information available so far to make the reservations at the next time slot. Again, since the reservation $A^t$ is selected before the job request $B^t$ is known, we will not be able to reach the optimal solution to  $(\ref{opt-pb-on})$, but instead, we opt to attain a total reservation cost that is not too large compared to some benchmark that knows the job requests in advance, and meanwhile, to ensure that the constraint is asymptotically satisfied on average over large time horizons. Specifically, we assume the following an "ideal" optimization problem where the values of job requests $\{B^t\}_{t=1}^T$ are known in advance for the entire horizon $T$: 
	\begin{equation} 
		\begin{split}
			\left\{
			\begin{tabular}{l}
				$\min\limits_{A\in\mathcal{R}} \sum\limits_{t=1}^T C(A)$,\\
				s.t.$\sum\limits_{s=1}^tC_0(A,B^s)\leq v$, $\forall t\leq T$.
			\end{tabular}
			\right.
		\end{split}
		\label{opt-pb-hind-no-rand}
	\end{equation}
	Denote by $A^*$ the optimal static solution to $(\ref{opt-pb-hind-no-rand})$. Therefore, the regret of not selecting $A^*$ over the horizon $T$ is defined by
	\begin{align}
		R_T= \sum_{t=1}^T \big( C(A^t)- C(A^*)\big). 
		\label{reg-rand}
	\end{align}
	We present in this paper a randomized algorithm that produces a sub-linear regret in time. This corresponds to a time average regret converging towards zero as the time reaches infinity. The intuition behind this is that the algorithm is learning and improves its performance as it accumulates knowledge. 
	\begin{rem}
		One can consider a different formulation of the problem by minimizing the total cost of reservation, violation, and transfer; see \cite{b24}. 
	\end{rem} 
	\begin{figure}
	\center
	\includegraphics{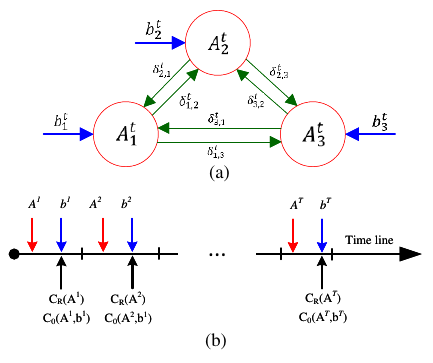}
	\caption{(a) The network model for three nodes at time $t$. (b) The reservation, job requests, and cost calculation timeline.}
	\label{fig:nodes}
	\end{figure}
	
	\section{Online randomized reservations} 
	\label{rand-formulation-sec}
	Notice that the reservation set $\mathcal{R}$ is discrete. Therefore, the classical gradient type algorithms (see e.g. \cite{b17}) might not be directly applicable since no convexity assumption is possible. To overcome this issue, we introduce randomization in the decision process where, at every time step, the network administrator chooses a probability distribution $P^t$ over the set $\mathcal{R}$ of reservations, based on the past job requests, and then draws the reservation $A^t$ randomly according to $P^t$. Especially, define 
	\begin{align*}   
		\mathcal{P}(\mathcal{R})=\bigg\{P:=(p_a)_{a\in\mathcal{R}}\in\mathbb{R}^{|\mathcal{R}|}_{+}:\mbox{ $\sum\limits_{a\in\mathcal{R}}p_{a}=1$}\bigg\},
	\end{align*} 
	as the space of probability distributions over $\mathcal{R}$. Moreover, define the expected reservation cost as $
	\mathbb{E}_{P^t}[C(A^t)]=\sum_{a\in\mathcal{R}}p^t_{a} C(a)$. Similarly, define the conditional expected transfer and violation costs given that the job request $B^t$ are known, respectively, by $\mathbb{E}_{P^t}[C_T(A^t, B^t)]=\sum_{a\in\mathcal{R}}p^t_{a}C_T(a, B^t)$ and $	\mathbb{E}_{P^t}[C_V(A^t, B^t)]=\sum_{a\in\mathcal{R}}p^t_{a}C_V(a, B^t)$. Since the reservations $\{A_t\}_{t=1}^T$  are generated randomly, the regret is now a random variable. Let us thus introduce the \textit{expected} regret $\tilde{R}_T$ defined in terms of expected costs. More precisely, conditioned on the fact that the values of job request $\{B^t\}_{t=1}^T$ are known in advance, let the fixed probability distribution $P^*=\argmin\limits_{P\in\mathcal{P}(\mathcal{R})}\sum\limits_{t=1}^T \mathbb{E}_{P}[C(A^t)]$ such that $\sum\limits_{s=1}^t\mathbb{E}_{P}[C_0(A,B^s)-v]\leq 0,\forall t\leq T$ (in analogy with the fixed allocation $A^*$ solution of $(\ref{opt-pb-hind-no-rand})$). Then, one defines the expected \textit{regret} of not playing according to $P^*$ over the horizon $T$ as $\tilde{R}_T= \sum\limits_{t=1}^T \big( \mathbb{E}_{P^t}[C(A^t)] - \mathbb{E}_{P^*}[C(A^t)]\big)$.

	\section{Exponentially weighted algorithm with long-time constraints}   
	\label{algo-sec}
	A classical approach for solving constrained convex optimization problems is the Lagrange multipliers method which consists of finding a saddle point of the Lagrangian function (see, e.g. \cite{b3}). Online versions of these methods have been proposed recently (see, e.g., \cite{b8}).
	
	We propose a different approach by replacing the minimization step in the saddle-point methods with an exponentially weighted step. Our idea is inspired by the exponentially weighted average algorithm introduced in \cite{b10} and \cite{b11}, which assigns a weight to each reservation vector based on its past performances. Namely, one wants to perform as well as the best reservation vector. Thus, one weights each reservation vector with a non-increasing function of its past cumulative costs, (\textit{cf} \cite[Chap. 4.2]{b22}). The novelty of the present paper is that it extends the exponentially weighted algorithm to incorporate time-varying constraints. This is done by introducing a step term in the exponent that takes into account the past cumulative constraint violations. In particular, one generates the reservations randomly according to the sequence of probability distributions $\{P^t\}_{t\geq 0}$ constructed, at each $t=1,2,\ldots$
	\begin{equation}
		\begin{split}
			P^t_a=\frac{w^t(a)}{\sum_{a\in\mathcal{A}}w^t(a)}, 
		\end{split}
	\end{equation}
	where the weights are given as 
	\begin{equation}
		\begin{split}
			w^{t+1}(a)=\exp\bigg\{-\eta \sum_{s=1}^{t}\big(C(a)+\lambda \digamma(s,a)\big)\bigg\},
		\end{split}
	\end{equation}
	with
	\begin{align}
		\digamma(s,a)= \bigg[\frac{1}{s}\sum_{r=1}^s(C_0(a,B^r)-v ) \bigg]^+, 
		\label{digamma}
	\end{align}
	$\eta>0$ a positive parameter, and $B^1,\ldots,B^{t-1}$ are the values of job requests observed so far. Moreover, $\lambda>0$ can be seen as the analogous of a fixed \textit{Lagrange multiplier} that calibrates the weight to put on the constraints. Furthermore, we initialize the algorithm with a uniform distribution, i.e. $w^1(a)=1$ for all $a\in\mathcal{A}$. Thus, the following recursion is straightforward:  
	\begin{align}
		&w^{t+1}(a)=w^{t}(a) \exp\big\{-\eta \big(C(a)+\lambda\digamma(t,a) \big)\big\}. 
		\label{rec-exp}
	\end{align}
	The advantage of our approach is twofold. First, one considers the whole history rather than the latest arrival. Steeply, by doing so, one avoids solving minimization problems which can be computationally heavy. We summarize our approach in Algorithm \ref{sad-algo}. The selection of parameters $\eta$ and $\lambda$ is discussed in Section \ref{per-sec}. Notice that the time-average cumulative constraint violations function introduced in $(\ref{digamma})$ is similar to the one used in \cite{b21} to extend \textit{Follow The Regularized Leader} algorithm to encompass time-varying constraints. 
	\section{Performance analysis}
	\label{per-sec}
	We analyze in this section the performance of Algorithm \ref{sad-algo}, first in terms of regret, and then in terms of cumulative constraint violations. Define the feasible sets $\mathcal{G}_t=\big\{a\in\mathcal{A}:\sum_{r=1}^s(C_0(a,B^r))\leq vs, \forall s\leq t \big\}$, and $\mathcal{E}_t=\big\{P\in\mathcal{P}(\mathcal{A}):\sum_{r=1}^s\mathbb{E}_P[C_0(a,B^r)]\leq vs, \forall s\leq t \big\}$. Moreover, we assume the following in the sequel: 
	\begin{enumerate}
		\item Cost boundedness: there exists a constant $\Theta>0$ such that, for all $a,b\in\mathcal{R}$,
		\begin{align*}
			\big|C(a)\big|\leq \Theta,\big|C_T(a,b)\big|\leq \Theta,\mbox{ and } \big|C_V(a,b)\big|\leq \Theta. 
		\end{align*}
		\item The subspaces $\mathcal{G}_t$ and $\mathcal{E}_t$ are not empty for all $t\geq 1$.  
	\end{enumerate} 
	\begin{algorithm}[h]
		\caption{Exponentially weighted method for online reservation with time-varying constraints}
		\begin{itemize}
			\item Initialize the values of $P^0\in\mathcal{P}(\mathcal{R})$, $B^0\in\mathcal{R}$, and $\lambda\geq 0$ 
			\item At each step $t=1,2,\ldots$  
			\begin{itemize}
				\item Observe $B^{t-1}$. 
				\item For each $a\in\mathcal{A}$ compute:
			\end{itemize}
			\begin{equation*}
				\begin{split}
					w^{t}(a)=\exp\bigg\{-\eta \sum_{s=1}^{t-1}\big(C(a)+\lambda \digamma(s,a)\big)\bigg\},
				\end{split}
			\end{equation*}
			\begin{align*}
				P^t_a=\frac{w^t(a)}{\sum_{a\in\mathcal{A}}w^t(a)}, 
			\end{align*}
			\item Generate $a^t$ according to $P^t$ and pay the costs
		\end{itemize}   
		\label{sad-algo}
	\end{algorithm}
	
	\subsection{Regret bound} 
	Let $\{P^t\}_{t\geq 1}$ be the sequence of probability distribution produced by Algorithm \ref{sad-algo} and let $\{A^t\}_{t\geq 1}$ be the sequence of corresponding random reservations. Recall that the \textit{regret} over finite horizons $T>0$, given by the formula $(\ref{reg-rand})$, quantifies the difference between the cumulative cost incurred by the exponentially weighted algorithm and the cumulative cost that would have been experienced if the algorithm had possessed perfect information and could have made the best-fixed decisions in hindsight. We show next that the regret is \textit{sublinear} in $T$ with high probability, regardless of the value of the job request sequence. 
	\begin{theo}
		For any $0<\delta<1$, the regret related to Algorithm \ref{sad-algo} satisfies, with a probability at least $1-\delta$,
		\begin{equation}
			\begin{split}
				R_T &\leq   \kappa T+\frac{\log |\mathcal{A}|}{\eta}+\Theta\sqrt{\frac{1}{2} \log (\delta^{-1}) T },
			\end{split}
			\label{reg-bound}
		\end{equation} 
		where $\kappa=\frac{\eta}{8} (1+2 \lambda)^2\Theta^2$. In particular, if $\eta=\frac{1}{\sqrt{T}}$, then $R_T=\mathcal{O}(\sqrt{T})$ (sublinear in $T$). 
		\label{reg-theo}
	\end{theo}
	
	\proof First, notice that 
	\begin{align*}
		\log\frac{\sum_{a\in\mathcal{A}}w^{t+1}(a)}{\sum_{a\in\mathcal{A}}w^1(a)}&\geq-\eta \min_{a\in\mathcal{G}_t}\sum_{s=1}^t C(a)-\log |\mathcal{A}|. 
	\end{align*} 
	Moreover, 
	\begin{align}
		\log \frac{\sum_{a\in\mathcal{A}}w^{t+1}(a)}{\sum_{a\in\mathcal{A}}w^1(a)}=\sum_{s=1}^t\log\frac{\sum_{a\in\mathcal{A}}w^{s+1}(a)}{\sum_{a\in\mathcal{A}}w^{s}(a)}.
		\label{log-prop}
	\end{align}
	Furthermore, define the random variable $A^s$ such that $\mathbb{P}(A^s=a)=\frac{w^s(a)}{\sum_{a\in\mathcal{A}}w^{s}(a) }$. Thus, by $(\ref{rec-exp})$ one obtains
	\begin{align*}
		&\log\frac{\sum_{a\in\mathcal{A}}w^{s+1}(a)}{\sum_{a\in\mathcal{A}}w^{s}(a)}=\log\mathbb{E}_{P^s}\big[\exp\big\{-\eta \big(C(a)\\
		&\qquad\qquad\qquad\qquad\qquad\qquad\qquad+\lambda \digamma(s,a)\big) \big\} \big]. 
	\end{align*}
	Now, by the assumptions above one deduces that $C(a)+\lambda \digamma(s,a)\leq  (1+2 \lambda)\Theta$. Therefore, by Hoeffding's lemma (\cite[Lemma 2.6]{b31}), one gets 
	\begin{align*}
		\log\frac{\sum_{a\in\mathcal{A}}w^{s+1}(a)}{\sum_{a\in\mathcal{A}}w^{s}(a)}&\leq -\eta\mathbb{E}_{P^s}[C(a)+\lambda \digamma(s,a)]+\kappa.
	\end{align*} 
	Thus, 
	\begin{align*}
		-\eta \min_{a\in\mathcal{G}_t}\sum_{t=1}^T C(a)&\leq -\eta \sum_{t=1}^T \mathbb{E}_{P^t}[C(a)+\lambda \digamma(t,a)]\\
		&\quad+\kappa T+\log |\mathcal{A}|.
	\end{align*}
	Now, notice that $C(A^*)= \mathbb{E}_{P^*}[C(a)]$, from which we get 
	\begin{equation*}
		\begin{split}
			\sum_{t=1}^T \mathbb{E}_{P^t}[C(a)+\lambda \digamma(t,a)]&- \mathbb{E}_{P^*}[C(A)]\leq \kappa T+\frac{\log |\mathcal{A}|}{\eta}.
		\end{split}
	\end{equation*}
	Now, since $\digamma(t,a)\geq 0$ for all $t\geq 0$ one obtains that the \textit{expected} regret satisfies
	\begin{equation}
		\begin{split}
			\tilde{R}_T
			&\leq  \kappa T+\frac{\log |\mathcal{A}|}{\eta}.
			\label{exp-reg-bound}
		\end{split}
	\end{equation}
Suppose that the values of job requests are $\{B^t\}_{t=1}^T=\{b^t\}_{t=1}^T$. Moreover, let $A^*$ be the optimal solution to the following  constrained combinatorial optimization problem   
\begin{equation} 
	\begin{split}		\min\limits_{A\in\mathcal{G}_t}\sum\limits_{t=1}^T C_R(A),		
	\end{split}
\end{equation} 
Therefore, the regret of not playing $A^*$ over the horizon $T$ is given by
\begin{align*}
	R_T= \sum_{t=1}^T \big( C_R(A^t)- C(A^*)\big). 
\end{align*}
Define the random variables 
\begin{align*}
	X^t= C_R(A)-\mathbb{E}_{P^t}[C_R(A^t)]. 
\end{align*}
Then, by the assumption above, $|X^t|\leq \Theta$. Therefore, $\{X_t\}_{t\geq 1}$ is a sequence of bounded martingales differences, and $Y_T=\sum_{t=1}^TX_t$ is a martingale with respect to the filtration $\mathcal{F}_T=\sigma(A^t,t\leq T)$. Thus, by Hoeffding-Azuma's inequality \cite{b32} one gets that, for any $0<\delta<1$,    
\begin{align*}
	\mathbb{P}\bigg(Y_T\leq \sqrt{\frac{1}{2} \log (\delta^{-1}) T \Theta^2}\bigg)\geq 1-\delta.
\end{align*}  
Thus, with probability at least $1-\delta$, one has 
\begin{align*}
	\sum_{t=1}^T C_R(A^t) \leq \sum_{t=1}^T \mathbb{E}_{P^t} [C_R(A^t)]+\sqrt{   \frac{1}{2} \log (\delta^{-1}) T \Theta^2}. 
\end{align*} 

Now, define the probability distribution $P_*\in\mathcal{P}(\mathcal{A})$ as the optimal solution to the following optimization problem  
\begin{equation} 
	\begin{split}		\min\limits_{P\in\mathcal{E}_t}\sum\limits_{t=1}^T \mathbb{E}_P[C_R(A)]
	\end{split}
	\label{p-star-stat-K}
\end{equation}  
if the values of job requests $\{B^t\}_{t=1}^T=\{b^t\}_{t=1}^T$ were known is advance. Then, since $\delta_{A^*}\in\mathcal{E}_t$,  $\sum_{t=1}^T\mathbb{E}_{P^*}[C_R(A)]\leq \sum_{t=1}^T C_R(A^*)$. Therefore, 
\begin{equation}
	\begin{split}
		R_T &\leq \sum_{s=1}^t \mathbb{E}_{P^s}[C_R(A^s)]- \inf_{P\in\mathcal{E}_t}\sum_{s=1}^t \mathbb{E}_P[C_R(A)]\\
		& \qquad\qquad+\Theta\sqrt{\frac{1}{2} \log (\delta^{-1}) T }. 
	\end{split}
\end{equation}
Using $(\ref{exp-reg-bound})$ gives (\ref{reg-bound}). \carre

	\subsection{Violation constraint bound}
	We analyze next the expected cumulative time average of the constraint violations.
	\begin{theo}
		The sequence $\{P^t\}_{t=1}^T$ produced by Algorithm \ref{sad-algo} satisfies: 
		\begin{equation}
			\begin{split}
				\frac{1}{T}\sum_{t=1}^T\mathbb{E}_{P^t}[\digamma(t,a)] &\leq \frac{1}{\lambda }\min_{a\in\mathcal{G}_t} C(a)+\frac{\log |\mathcal{A}|}{T\lambda\eta}+\frac{1}{\lambda} \kappa. 
				\label{const-upp-bound}
			\end{split}
		\end{equation}
		In particular, if $\eta=\frac{1}{\sqrt{T}}$, $\lambda=o(T^{1/4})$ then 
		\begin{align}
			\lim_{T\rightarrow\infty}\frac{1}{T}\sum_{t=1}^T\mathbb{E}_{P^t}[\digamma(t,a)]\rightarrow 0. 
			\label{lim-const}
		\end{align}
	\end{theo}
	\proof  First, notice that 
	\begin{align*}
		\log\frac{\sum_{a\in\mathcal{A}}w^{t+1}(a)}{\sum_{a\in\mathcal{A}}w^1(a)}&\geq -\eta \min_{a\in\mathcal{A}}\sum_{s=1}^t (C(a)+\lambda\digamma(s,a))\\
		&\quad-\log |\mathcal{A}|.
	\end{align*} 
	On the other hand, by $(\ref{log-prop})$ and $(\ref{rec-exp})$, one obtains
	\begin{align*}
		\log\frac{\sum_{a\in\mathcal{A}}w^{s+1}(a)}{\sum_{a\in\mathcal{A}}w^{s}(a)}&=\log\mathbb{E}_{P^s}\big[\exp\big\{-\eta (C(a)\\
		&\qquad\qquad\qquad+\lambda\digamma(s,a)) \big\} \big]. 
	\end{align*}
	Therefore, by Hoeffding's lemma (see, e.g. \cite[Lemma 2.6]{b31}) and the above assumptions, one gets 
	\begin{align*}
		\log\frac{\sum_{a\in\mathcal{A}}w^{s+1}(a)}{\sum_{a\in\mathcal{A}}w^{s}(a)}&\leq -\eta\mathbb{E}_{P^s}[(C(a)+\lambda\digamma(s,a))]+\kappa.
	\end{align*} 
	One thus deduces $(\ref{const-upp-bound})$.
	Finally, taking $\eta=\frac{1}{\sqrt{T}}$ and $\lambda=o(T^{1/4})$ leads to $(\ref{lim-const})$. \carre
	
	\section{Numerical experiments and comparison with reinforcement learning} 
	\label{num-sec} 
	To show the performance of our algorithm, we compare it with a reinforcement learning approach tailored to the network resource allocation problem. Then, we test the two algorithms over two scenarios. The first consists of piecewise-fixed input requests and the next step is randomly generated inputs.
	
	\subsection{Reinforcement learning algorithm}
	\label{sec:rl}
	Reinforcement learning (RL) is an approach to finding a control policy through an \textit{educated} interaction of the system with the environment (\textit{cf} \cite{b33}). In RL, an agent interacts with an environment by making decisions and receiving a reward. We model the environment as an $n$-armed bandit, where the agent has to find the option with the highest pay-off among its many options, coupled with a fuzzy inference system (FIS). The input of the decision maker is then the average of input requests to each node in the past $T_H$ time steps, and the output is the number of reserved resources. The decision maker has a centralized structure. The output of the FIS is a scalar, which then will be mapped to a vector, in which each element represents the reservation for each node. We use a Takagi-Sugeno fuzzy structure, with discrete output as in \cite{b34}. In addition, we use a training approach similar to the fuzzy actor-critic learning (FACL) algorithm as in \cite{b34}. We modify the FACL algorithm in a way that it can solve the $n$-armed bandit problem. The RL approach we propose is shown in Algorithm \ref{al:rl}.

	\subsection{Piecewise-Fixed Requests}
	\label{ssec:fixed_requests}
	First, as our benchmark for testing, we take a network with three nodes, and all the nodes are connected. To put our written code to the test, we send a piecewise fixed request to each node. The input requests will be constant for the first 500 steps. At $t=1000$ steps, the input requests will be changed again. From $t=1001$ steps to $t=5000$ steps the inputs will be fixed. The reservation cost for each node is $f_R=0.05x^2$, where $x$ is the reservation. The violation cost for each node is $f_V=0.05x^2$, where $x$ is the number of requests, more than a node's reservation. The transfer cost is $f_T^{nm}=k_{nm}x^2$, which is the transfer cost of $x$ jobs from node $n$ to node $m$. The coefficients $k_{nm}$ are shown in Table \ref{tab:k}
	\begin{table}[htbp]
		\centering
		\caption{The transfer cost coefficients ($k_{nm}$)}
		\begin{tabular}{cccc}
			\hline
			$k_{nm}$ & \multicolumn{3}{c}{$n$} \bigstrut\\
			\hline
			\multirow{3}[2]{*}{$m$} & 0     & 0.02  & 0.03 \bigstrut[t]\\
			& 0.02  & 0     & 0.02 \\
			& 0.03  & 0.02  & 0 \bigstrut[b]\\
			\hline
		\end{tabular}%
		\label{tab:k}%
	\end{table}%
	In Algorithm \ref{al:rl}, $\alpha$ is the learning rate and is set to $0.995$ and $\varepsilon$ is 0.5. In addition, each actor's input has 11 triangular membership functions. The term $T_H$ is set to $100$ steps and $T_B$ is set to $50$ iterations. To demonstrate the performance we used the time average regret $\frac{1}{T}R_T$. It should be mentioned that the time horizon $T$ is not known for the user. Thus, nature's optimization in \eqref{opt-pb-hind-no-rand} must be done for any possible $T$. To do the optimization for finding $A^*$, we use a genetic algorithm by using the built-in Matlab function \textit{ga}. In the simulation, the term $\lambda$ in \eqref{rec-exp} and the reward function in Algorithm 2 is set to 32. The performances of Algorithm \ref{sad-algo} and Algorithm \ref{al:rl} are shown in Fig. \ref{fig:constant_request}.
	\begin{figure}
		\centering
		\includegraphics{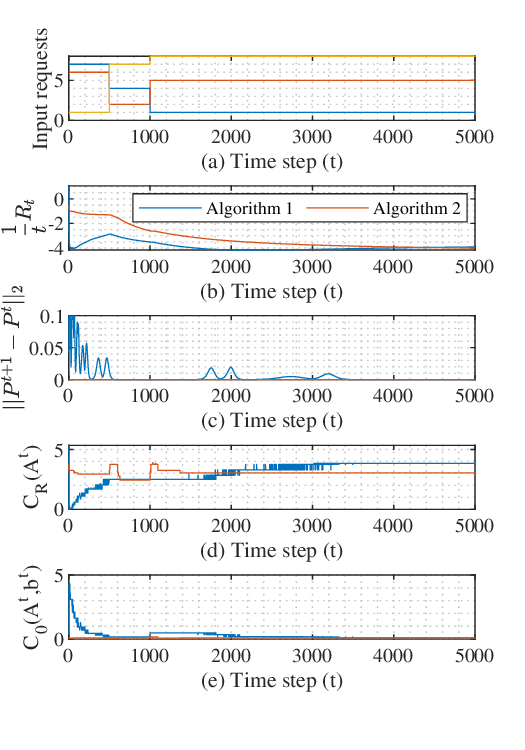}
		\caption{(a) The input request to each node. (b) The average regret of both methods. (b) The Euclidean distance between two consecutive probability distributions produced by Algorithm \ref{sad-algo}. (c) Reservation cost ($C(A^t)$) for each time slot $t$. (d) The blocking cost ($C_0(A^t,B^t)$)}.
		\label{fig:constant_request}
	\end{figure}
	Fig. \ref{fig:constant_request} (a) shows the input requests to each node. Fig. \ref{fig:constant_request} (b) compares the average regret for the RL method (Algorithm \ref{al:rl}) and the exponentially weighted method (Algorithm \ref{sad-algo}). It shows for most of the time horizon, the average regret of Algorithm 1 is less than Algorithm 2. This is evidence that Algorithm \ref{sad-algo} outperforms Algorithm \ref{al:rl}. Fig. \ref{fig:constant_request} (c) shows the Euclidean distance between consecutive probability distributions for the exponentially weighted method. Notice that the RL method does not calculate $P^t$. Fig. \ref{fig:constant_request} (d) shows the reservation cost. It is observed that Algorithm 2 responds fast to each change in the environment. The reason is that the RL method only sees the past $T_H$ time steps, while Algorithm 1 sees the entire history. Finally, Fig. \ref{fig:constant_request} (e) shows the blocking cost. One can observe that the RL method can keep the blocking cost at $v$ more effectively than the exponentially weighted method. For instance, after the change at $t=1000$ steps, the exponentially weighted method's blocking cost goes up to $0.49$ and returns to $0.1$ after $2000$ steps, while the RL method's blocking cost jumps to 0.17 and returns to $0.1$ after $100$ steps.
	\begin{algorithm}[h]
		\caption{Discrete Fuzzy Actor-Critic Learning Algorithm \cite{b34} for $n$-Armed Bandit Problem}
		\begin{itemize}
			\item Initialize the \textit{actor}, and then the \textit{critic} with Gaussian distribution $\mathcal{N}(0,1)$.
			\item Set $\alpha$: the \textit{learning rate}, and $\varepsilon$ the \textit{exploration-exploitation rate}.
			\item Set the \textit{back processing iterations} ($T_B$), and the \textit{time history} ($T_H$).
			\item At each step $t=1,2,\ldots$, observe $B^{t-1}$ and  
			\begin{itemize}
				\item Calculate the FIS inputs:\\ 
				$M_i(t,T_H)=Mean(B^{max(1,t-T_H)}_i,\ldots,B^{t}_i), \\ \forall i=1,\ldots,G \text{(Number of nodes)}$.
				\item The actor's output parameters ($\omega$) are selected. The parameters corresponding to the critic's largest action given a rule are selected.
				\item $a^t=FLC((M_1(t,T_H),\ldots,M_G(t,T_H)),\omega)$.
				\item Convert $a^t$ to $A^t$ (Since $a^t$ is a scalar and must be mapped to a vector).
				\item Do the back processing for $T_B$ times.
				\begin{itemize}
					\item for $j=1,\ldots,T_B$
					\item Actor's output parameters ($\omega$) are selected. The parameters corresponding to the critic's largest action given a rule are selected. After selection, $\varepsilon\times100\%$ of the input parameters are randomly changed.
					\item $a^t=FLC((M_1(t,T_H),\ldots,M_G(t,T_H)),\omega)$.
					\item Convert $a^t$ to $A^t$.
					\item Calculate Reward:\\ $R(A^t,B^{t-1})=-C(A^{t-1})-\lambda C_0(A^t,B^{t-1})$.
					\item Update critic by having $R$ and $\alpha$.
				\end{itemize}
				\item Calculate the cost $C(A^t,B^t)$.
			\end{itemize}
		\end{itemize}   
		\label{al:rl}
	\end{algorithm}
	\subsection{Randomly generated requests}
	Now, we repeat the previous experiment for a sequence of randomly generated requests. The total time horizon is set to $T=5000$ steps divided into 22 regions based on a Poisson process with a mean of $250$ time steps. In each region, the requests are randomly generated by a Poisson distribution. Each node in each region has a unique integer as the mean that is selected from \{2,3,4,5,6\}. If the generated request is more than $10$, it will saturated to $10$. We then apply Algorithm \ref{sad-algo} and \ref{al:rl} and compare their performances by observing their regrets.
	\begin{figure}
		\centering
		\includegraphics{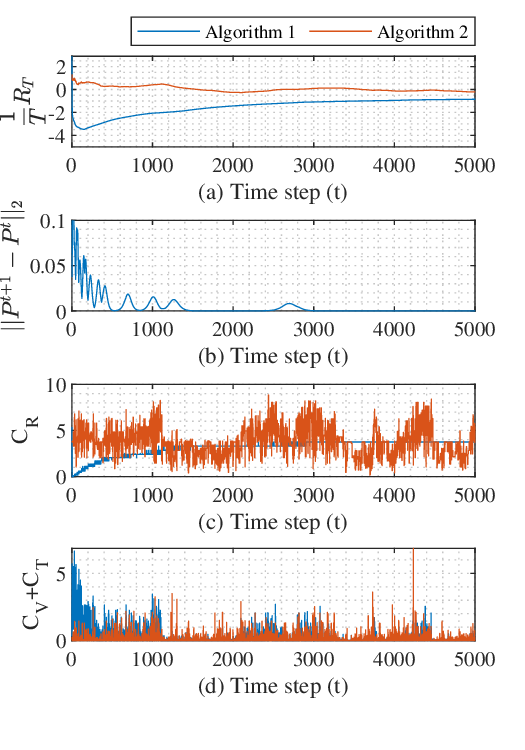}
		\caption{(a) The input request to node 1. The rest inputs are analogous. (b) The average regret of both methods. (b) The Euclidean distance between two consecutive probability distributions produced by Algorithm \ref{sad-algo}. (c) Reservation cost ($C(A^t)$) for each time slot $t$. (d) The blocking cost ($C_0(A^t,B^t)$)}.
		\label{fig:random_request}
	\end{figure}
	Fig. \ref{fig:random_request} (a) shows the input requests to node 1. The input requests to node 2 and 3 are analogous. Fig. \ref{fig:random_request} (b) compares the performances of two algorithms. It is shown that for the whole time horizon Algorithm 1 has less regret than Algorithm 2. Moreover, Fig. \ref{fig:constant_request} (c) shows the Euclidean distance between consecutive probability distributions for the exponentially weighted method. Furthermore, Fig. \ref{fig:constant_request} (d) demonstrates an interesting fact about the difference between the two approaches where the exponentially weighted method reservation cost is almost the same throughout the time horizon, whereas the RL method often tries to adapt and comes with a new reservation. The exponentially weighted method reaches an atomic probability density function. Again, this is evident that Algorithm \ref{sad-algo} outperforms Algorithm \ref{al:rl}. Finally, Fig. \ref{fig:constant_request} (d) shows the blocking cost of the two approaches. It shows that except for the time interval of $0$ to $1000$ steps, the blocking cost is almost the same for the two algorithms. In the initial 1000 time steps, Algorithm 1 has not converged yet.
	
	\section{Conclusion}
	In this paper a specific online resource allocation problem with dynamic job transfers. We proposed a new randomized online algorithm based on an exponentially weighted method that encompasses time-varying constraints. We showed that our method has a sublinear in time regret expected constraint violations. This indicates that it adapts and improves its performance as it learns from the past. Moreover, we compared our new approach with a reinforcement learning algorithm where we observed that our proposed method outperforms the latter when the job requests sequence is varying. Indeed, in both tested data streams, the regret is lower for the exponentially weighted method presented in this paper, while the average constraint violations decrease over time.

		\section{Appendix}
		
		\begin{figure}
			\centering
			\includegraphics{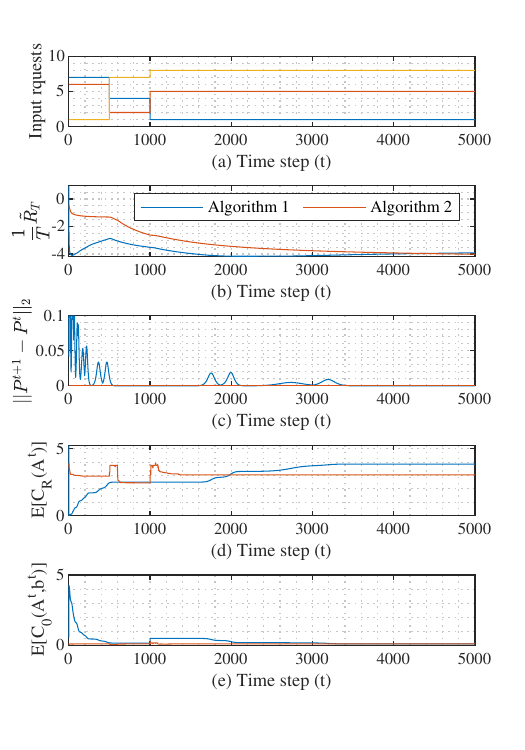}
			\caption{(a) The regret of both methods. (b) The Euclidean distance between two consecutive probability distributions produced by Algorithm \ref{sad-algo}, (c) the corresponding reservation const ($E[C(A^t,B^t)]$), and (d) the blocking cost ($E[C_0(A^t,B^t)]$)}.
			\label{fig:constant_request}
		\end{figure}
		
		\begin{figure}
			\centering
			\includegraphics{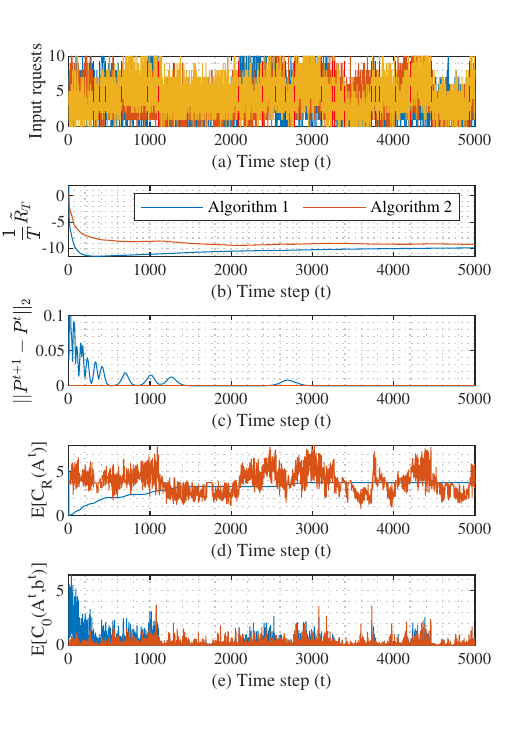}
			\caption{(a) The regret of both methods. (b) The Euclidean distance between two consecutive probability distributions produced by Algorithm \ref{sad-algo}, (c) the corresponding reservation const ($E[C(A^t,B^t)]$), and (d) the blocking cost ($E[C_0(A^t,B^t)]$)}.
			\label{fig:constant_request}
		\end{figure}

\end{document}